\input amstex
\input epsf
\define\W{\Cal{W}}
\define\ltilde{\sim}

\define\tna{T_n^A}

\define\tnb{T_n^B}

\define\tns{T_n^S}
\define\ncbn{\operatorname{NC}_n^B}
\define\ncns{\operatorname{NC}_n^S}

\define\mi{M_n^{(i)}}
\define\mii{M_n^{(ii)}}

\define\upa{\uparrow\!}
\define\downa{\downarrow\!\!}

\define\ncbd{\operatorname{NC}_n^S}

\define\ncna{\operatorname{NC}_n^A}
\define\ncnb{\operatorname{NC}_n^B}

\NoBlackBoxes
\documentstyle{amsppt}

\topmatter
\title Tamari lattices and noncrossing partitions 
in type $B$ and beyond \endtitle
\rightheadtext {Tamari lattices and noncrossing partitions}
\author Hugh Thomas \endauthor

\abstract 
The usual, or type $A_n$, Tamari lattice is a partial order on $T_n^A$,
the triangulations of an $(n+3)$-gon.
We define a partial order on $T^B_{n}$, the set of 
centrally symmetric triangulations of a $(2n+2)$-gon. We show that
it is a lattice, 
and that it shares certain other nice properties
of the $A_n$ Tamari lattice, and therefore that it deserves to be 
considered the $B_n$ Tamari lattice.  

We define a bijection between
$T_n^B$ and the noncrossing partitions of type $B_n$ defined by Reiner.    
For $S$ any subset of $[n]$, Reiner defined a pseudo-type
$BD^S_n$, to which is associated a subset of the noncrossing partitions
of type $B_n$.  
We 
show that the elements of $T_n^B$ which correspond to the noncrossing
partitions of type $BD^S_n$ posess a lattice structure induced from their
inclusion in $T_n^B$.
\endabstract

\endtopmatter

\document

\head Introduction \endhead

The usual (or type $A_n$) Tamari lattice is a partial order on 
$\tna$,  the set of triangulations of an $(n+3)$-gon.  
By a triangulation of a polygon, we mean a division of the polygon 
into
triangles by connecting pairs of its vertices with straight lines which
do not cross in the interior of the polygon.  
The purpose of this paper is to define and investigate the properties
of an analogous lattice defined on centrally symmetric triangulations of
a $2n+2$-gon, which we call the type $B_n$ Tamari lattice.  More explanation
for why these lattices should be associated to the reflection groups of
types $A_n$ and $B_n$ 
will be given below.  

We begin by reviewing some features of the type $A_n$ Tamari lattice.  
Conventionally,
we will number the vertices of our $(n+3)$-gon        
clockwise 
from 0 to $n+2$, with a long top edge connecting vertices
$0$ and $n+2$.  An example triangulation is shown in Figure 1 below.

Let $S \in \tna$.  As in [Lee], we colour the chords of $S$ red and green,
as follows.  A chord $C$ of $S$ is the diagonal of a quadrilateral 
$Q(C)$ in $S$.
If $C$ is the diagonal of $Q(C)$
which
is connected to the vertex with the largest label, 
we colour it green; otherwise we colour it
red.  In Figure 1, the red chords are indicated by thick lines.  

$$\epsfbox{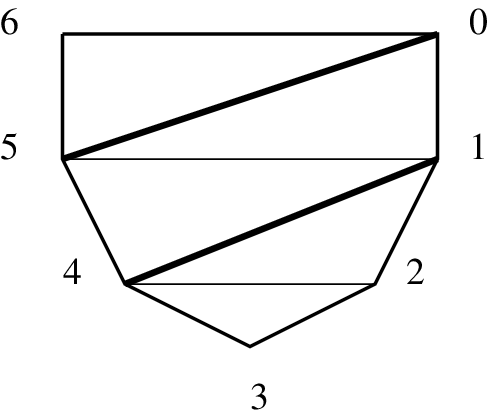}$$
$$\text{Figure 1}$$

We partially order $\tna$ by giving covering relations: $T$ covers $S$ if
they coincide except that some green chord in $S$ has been replaced by 
the other diagonal of $Q(C)$ (which is red).  This is one way to construct the 
the Tamari lattice, which was introduced in [Tam] and which has 
since been studied by several authors (see [HT, Pal, Mar, BW2]).  

Although this is not clear from the elementary description given here,
the Tamari lattice should be thought of as belonging to type $A$.  
One indication of why can be found in [BW2], where it is shown that
$\tna$ is a quotient of the weak order on the symmetric group $S_{n+1}$ (the
type $A_n$ reflection group). 
Another reason 
is that the elements of $\tna$ index clusters in the $A_n$
root system (see Fomin and Zelevinsky [FZ]).   
Once one has
the idea that the Tamari lattice belongs to type $A$,
it is natural to ask whether there exist Tamari lattices
in other types.  

For reasons which we shall go into further below, the $B_n$ triangulations,
denoted $\tnb$, 
are the triangulations of a centrally symmetric $(2n+2)$-gon which are 
themselves centrally symmetric (i.e.\ fixed under
a half-turn rotation).  
These triangulations have already appeared in the
work of Simion [Sim], and in [FZ]
where they index the clusters in the $B_n$ root system.  One goal of
our paper is to define a partial order on $\tnb$ and to prove that it is
a lattice.  The definition is analogous to that already given for the 
$A_n$ Tamari lattice: it is given in terms of covering relations, and
$S$ covers $T$ in $\tnb$ if $S$ is obtained from $T$ by replacing a symmetric
pair of chords $C, \bar C$ by the other diagonals of $Q(C), Q(\bar C)$.  
The details of the definition are a trifle complicated, so we defer them 
for the main body of the paper.  This definition was arrived at independently
and
more or less
simultaneously by Reading [Rea].  He has also proved that $\tnb$ is a lattice,
using a rather different approach.  Two alternative partial orders on
$\tnb$ with similar (but somewhat easier to describe) covering relations 
were suggested by Simion [Sim]; one is studied further in [Sa1].  
Since neither of these is a lattice, neither is completely satisfying
as a type $B$ analogue of the usual Tamari lattice.  

An important property of $\tna$ is that there is a natural (though not
order-preserving) bijection from $\tna$ to the (classical) noncrossing
partitions, $\ncna$.  The type $B_n$ version of noncrossing partitions,
$\ncnb$
was introduced by Reiner [Rei].  We show that there is a bijection
from $\tnb$ to $\ncnb$ similar to that from $\tna$ to $\ncna$. 

In the same paper where Reiner defined the noncrossing partitions of type
$B_n$, he defined a more general notion of noncrossing partitions of 
(pseudo-)type
$BD_n^S$, for $S\subset [n]$.  
The noncrossing partitions of type $BD_n^S$ are a subset of $\ncnb$.  Using
the above bijection between $\tnb$ and $\ncnb$, for any $S$, 
we describe the subset of
$\tnb$ which corresponds to the noncrossing partitions of type $BD_n^S$, and
we show that the order induced on this subset of $\tnb$ 
by its inclusion
in $\tnb$ gives it a lattice structure, which we call the 
Tamari lattice of pseudo-type $BD^S_n$.

We show that $\tnb$ (and the Tamari lattices of type $BD^S_n$ also) 
have an unrefinable chain of left modular
elements, a property also shared by the usual Tamari lattice [BS].  One 
consequence of this, due to Liu [Liu], is that these lattices have
EL-labellings.  Using these labellings, we show that, as for the usual 
Tamari lattice (see [BW2]), 
the order complex of any interval is
either homotopic to a sphere or contractible.  (This result on homotopy
types of order complexes of intervals was also
obtained by Reading [Rea].)

From the results in this paper one could proceed in two directions.  One 
direction is to consider the existence of Tamari lattices in all Coxeter
types.  The other direction is to investigate further the lattices defined
here, to see how many more of the properties of the usual Tamari lattice carry
over.

\head Type $B$ Triangulations \endhead

Recall that the $B_n$ Weyl group consists of signed permutations of 
$n$.  We can think of these as permutations of $\{1,\dots,n,\bar 1,
\dots,\bar n\}$  fixed under interchanging 
$i$ and $\bar i$ for all $1\leq i \leq n$.  By analogy,
 $B_n$ triangulations, $\tnb$,  are defined to
be type $A$ triangulations of a $(2n+2)$-gon
fixed under a half-turn.  
There is general consensus that this is the correct choice
of $B_n$ triangulation: see [Sim], [FZ].

We number the vertices of our standard $(2n+2)$-gon 
clockwise from $1$ to $n+1$ and then from $\overline{1}$
to $\overline{n+1}$.  A typical triangulation is shown in Figure 2.

We will frequently distinguish
two types of chords: {\it pure} and {\it mixed}.  A chord is pure if it
connects two barred vertices or two unbarred vertices; otherwise it is
mixed.  
For $S \in \tnb$, consider a  chord $C$ of $S$.  
The chord $C$ is the 
diagonal of a quadrilateral, which we denote $Q(C)$.
If $C$ is pure,
then we colour it red if $Q(C)$
contains another vertex of the same type as those of $C$ 
whose label is higher, 
and green otherwise.  If it is mixed, we colour it
red if $Q(C)$ contains an unbarred vertex whose label
is higher than the label of the unbarred vertex of $C$, 
or a barred vertex whose
label is higher than the 
label of the barred vertex of $C$.  Otherwise we colour it green. In
Figure 2, the red chords are indicated by thick lines.

$$\epsfbox{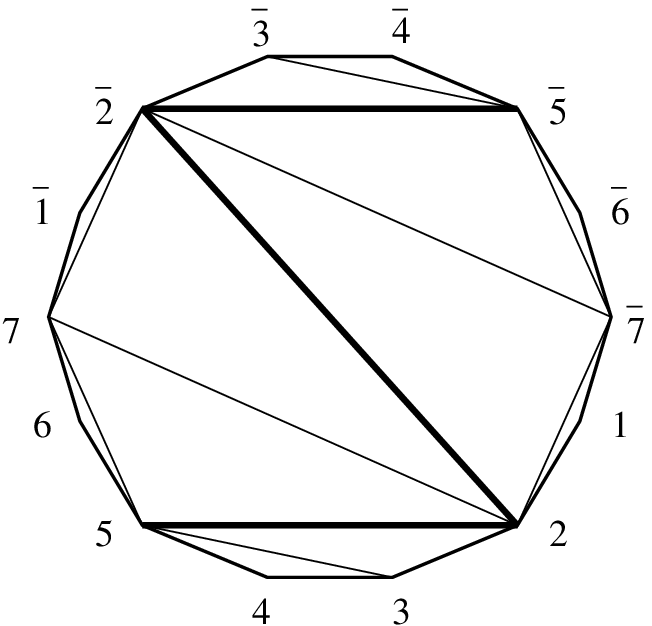}$$
$$\text{Figure 2}$$

For $C$ a chord, we write $\bar C$ for its symmetric partner (that is to say,
the image of $C$ under a half turn).  Observe that $C$ and $\bar C$ are
assigned the same colour.  

\proclaim{Lemma 1} Consider a chord $C$ in a triangulation $S$.  Let 
$S'$ be the triangulation obtained by replacing $C$
by $C'$, the other diagonal
of $Q(C)$, and also replacing $\bar C$ by
$\bar C'$.  Then the colours of $C$ in $S$ and $C'$ in $S'$ are opposite.
\endproclaim

\demo{Proof} The proof is just a case-by-case check of the possible 
configurations of the four vertices of $Q(C)$: all of one type,
three of one type and one of the other, 
or two of each type.
\qed \enddemo

We can now state the first main theorem of this paper (which, as was 
already remarked, was arrived at and proved independently and 
more or less simultaneously by
Reading [Rea]).  

\proclaim{Theorem 1} There is a lattice structure
 on $\tnb$ whose 
covering relations are given by $S \lessdot T$ iff 
$S$ and $T$ differ in that 
 green chords $C$, $\bar C$ in $S$ are replaced in $T$ by the other diagonals 
of
$Q(C)$ and $Q(\bar C)$ (which will be red).  Note that we allow 
$C=\bar C$ (i.e. $C$
being a diameter).  
We call this lattice the
$B_n$ Tamari lattice.  
\endproclaim
 
\demo{Proof} 
The proof of this theorem will take the rest of this section and all of the
next.  
We begin with a quick outline.  
We will associate to every triangulation $S\in \tnb$ a {\it bracket vector} 
$r(S)$ which is an $n$-tuple of elements from  $[0,n-1]\cup\{\infty\}$.  
We will define a partial order on $\tnb$ in terms of bracket vectors, and
then show
that its covering relations are as given in the statement of the theorem,
and that it is a lattice.

The first ingredient in our proof of Theorem 1 is some further 
analysis of the red and green chords of triangulations.  
Fix a triangulation $S$.  For $1\leq i \leq n$, look for a vertex of 
the polygon which is connected to $i$, starting at $\bar 1$ and searching
clockwise.  If none is found before reaching the vertex next counterclockwise
from $i$ (i.e.\ the vertex $i-1$, 
unless $i=1$, in which case the vertex $\overline{n+1}$), then $C_i(S)$ is
the edge segment connecting $i$ and the next vertex counterclockwise.  
Otherwise, if a vertex connected to $i$ was found, then
$C_i(S)$ is the chord of $S$ connecting $i$ to that vertex.
  Let 
$R(S)$ be the set of the $C_i(S)$ which are chords rather than edge 
segments, together with their symmetric partners.  

\proclaim{Lemma 2} For any triangulation $S$, $R(S)$ consists of 
the red chords; the chords not in $R(S)$ are green.  \endproclaim

\demo{Proof} Pick a chord in $R(S)$.  Since the colouring is symmetric, we
may assume that the chord 
is $C_i(S)$ for some $i$.  It follows that $Q(C_i(S))$ 
contains a vertex greater than $i$, and hence that 
$C_i(S)$ is red.  

Now consider a chord $C$ not in $R(S)$.  Suppose first that $C$ is pure;
we may assume that it connects $i$ and $j$ with $i>j$.  Since $C \ne C_i(S)$,
$C_i(S)$ divides $Q(C)$ from all those vertices with unbarred
labels greater than $i$, so $C$ is green.  Next suppose that $C$ is 
mixed, connecting $i$ and $\bar j$.  Now $C_i(S)$ divides $Q(C)$ 
from those vertices with unbarred labels greater than $i$,
and $\overline {C_j(S)}$ divides it from those vertices with barred labels greater
than $\bar j$.  So $C$ is green.  \qed\enddemo

\proclaim{Lemma 3} Let $M$ be a region of the $(2n+2)$-gon, that is to say,
the convex hull of some subset of the vertices of the $(2n+2)$-gon.
Then
there is a unique way to triangulate $M$ using only green chords.\endproclaim

\demo{Proof} The way to do it is as follows: connect every unbarred vertex
to the largest unbarred vertex, every barred vertex to the largest barred
vertex, and, if both exist, connect the largest barred and largest
unbarred vertices.  It is easy to see that all these chords are green.

Uniqueness is clear in the case where there is one type of vertex (that
is to say, barred or unbarred)
which appears at most once in the region.  So suppose 
we are not in this case, and fix a triangulation of $M$ using green
chords.  We wish to show that it is the triangulation defined in the
previous paragraph.  

Write $x$ for the largest unbarred vertex in $M$. 
Observe that there must be at least one mixed chord
in the triangulation, but no $i<x$ can have a mixed chord inside $M$ attached
to it, because this would imply by Lemma 2 that there was a red chord inside
$M$, contradicting our assumption.  Thus there is a mixed chord in the 
interior of $M$ which is connected to $x$, say $x\bar j$.

Now consider the triangle containing $x \bar j$, which is on the side 
of $x\bar j$ with the smaller unbarred
labels and the larger barred labels.  Suppose first that its third vertex
is barred.   
In this case, $x\bar j$ is not green, contradicting our
assumption.  So the third vertex must be unbarred, say $z$.  
Now $z\bar j$ cannot
be green.  Thus, it must be a boundary of $M$.  This implies that
$\bar j$ must be the greatest barred vertex of $M$.  Thus, we have shown
that our triangulation of $M$ contains the chord connecting the greatest
barred vertex and the greatest unbarred vertex.  This chord divides $M$
into two subregions which fall into the simple type (no more than one
barred vertex or no more than one unbarred vertex) for which uniqueness
is clear.  This establishes that the triangulation with which we began
must coincide with that described in the first paragraph of this proof. 
\qed \enddemo

The type $A$ analogue of Lemma 3 was proved in [Lee].  

\proclaim{Lemma 4} For any triangulation $S$, 
$S$ is the unique triangulation whose set of red chords
is exactly $R(S)$.  \endproclaim

\demo{Proof} Let $T$ be a triangulation whose set of red chords is $R(S)$.
The chords of $R(S)$ divide the $(2n+2)$-gon up into regions which are
triangulated by green chords of $T$, but by Lemma 3 there is a unique
way to do this, which must be that of $S$. So $T$ coincides with $S$.  
\qed\enddemo

\enddemo

\head  Bracket Vectors in types $A$ and $B$\endhead

We briefly recall some well-known facts about the type $A$ Tamari
lattice, which serve as motivation for our work in type $B$.  

Any triangulation $S \in \tna$ has a  bracket vector
$r(S)=(r_{1}(S), \dots, r_{n+1}(S))$.  
Let $v_i(S)$ be the least vertex attached to $i$.  Then 
$r_i(S)=i-1-v_i(S)$.  For example, the  bracket vector of the
triangulation shown in Figure 1 is (0,0,0,2,4).  
This approach to representing elements of the
Tamari lattice goes back to [HT], though we make some different choices
of convention here.  

\proclaim{Proposition 1 [HT]}
An $(n+1)$-tuple of positive integers is a  bracket vector for some
triangulation in $\tna$ iff it satisfies the following two properties: 

(i)  For $1\leq i<j\leq n+1$, $r_i\leq r_j-(j-i)$ provided 
$r_j-(j-i)$ is non-negative.
 
(ii) $0\leq r_i\leq i-1$.
\qed \endproclaim

The order relation on triangulations has a simple interpretation in terms
of  bracket vectors, which we summarize in the following proposition:

\proclaim{Proposition 2 [HT, Mar]} The lattice structure on $\tna$ can be 
described as follows:

(i)  $S\leq T$ iff 
$r_i(S)\leq r_i(T)$ for all $i$. 

(ii) $r_i(S\wedge T) =
\min(r_i(S),r_i(T))$.  

(iii) For $x$ any $n+1$-tuple of numbers satisfying only the second
condition of Proposition 1, there is a unique triangulation
$\upa\!(x)$
such that  for $S\in \tna$,  

$$r_i(S) \geq x_i \text{ for all $i$ iff } S\geq \upa\!( x).$$  

(iv)
$r(S\vee T)=\upa({\max}(r(S),r(T)))$, where max is taken coordinatewise.  
\qed\endproclaim

We now proceed to describe a similar construction in type $B$.
  To a triangulation $S \in \tnb$ we associate a  bracket vector
$r(S)=(r_1(S),\dots,r_{n}(S))$, as follows.
For $1\leq i \leq n$, 
let $v_i(S)$ denote the end other than $i$ of $C_i(S)$.  
If the counter-clockwise distance from $i-1$ to $v_i(S)$ is less than or 
equal to $n-1$, set $r_i(S)$ to be that distance.  Otherwise, set $r_i(S)=\infty$. 
Thus, the triangulation shown in Figure 2 has  bracket vector
$(0,\infty,0,0,2,0)$.  

\proclaim{Conventions regarding $\infty$} 
$\infty$ is considered to be greater than any integer.  $\infty$ plus 
an integer (or $\infty$) equals $\infty$.  
\endproclaim

\proclaim{Lemma 5} From the bracket vector $r(S)$, $C_i(S)$ can be determined
as follows:

(i) If $r_i(S)=0$ then $C_i(S)$ is the edge segment extending counter-clockwise
from $i$.

(ii) If $0<r_i(S)<n$ then $C_i(S)$ connects $i$ to the vertex $r_i(S)+1$ vertices
counter-clockwise from $i$.  

(iii) If $r_i(S)=\infty$ then $C_i(S)=i\bar j$ where $j$ is the least vertex satisfying
$r_j(S)-j\geq n-i$.  
\endproclaim

\demo{Proof} The cases $r_i(S)=0$ and $0<r_i(S)<n$ are obvious.

Suppose $r_i(S)=\infty$.  In this case,
certainly $C_i(S)=i\bar j$ for some $j$.  
It follows that $j \bar i$ is also a chord of
$S$.  
Thus $r_j(S)-j\geq n-i$.  
Now, for $m<j$, it cannot be that
vertex $m$ is connected to $\bar i$, since then $C_i(S)$ would be 
$i \bar m$, not $i \bar j$.  Also, it cannot be that $m$ is connected to any
vertex less than $\bar i$, since then $C_m(S)$ would divide $j$ from
$\bar i$.  Thus, $r_m(S)-m< n-i$.  It follows that $j$ is the least number
satisfying $r_j-j\geq n-i$, as desired.  
\qed
\enddemo

\proclaim{Corollary} 
The map from $\tnb$ to  bracket vectors is
injective. \endproclaim

\demo{Proof} Given $r(S)$, we can determine $C_i(S)$ for all $i$.  
Their union together with their symmetric partners gives $R(S)$, and,
by Lemma 4, determines $S$. \qed\enddemo

\proclaim{Proposition 3}  
$B_n$ 
 bracket vectors are $n$-tuples of symbols from $[0,n-1] \cup \{\infty\}$
characterized by the following two properties:

(i) For $1 \leq i<j\leq n$,
$r_i \leq r_j-(j-i) $ if $r_j-(j-i)$ is non-negative.

(ii) If $\infty>r_i\geq i$, then $r_{n+i-r_i}=\infty$.  

\endproclaim

\demo{Proof} Clearly a  bracket vector satisfies condition (i).  
Suppose that $\infty>r_i\geq i$.  Then $i$ is connected to $\overline{n+i-r_i}$.
By symmetry, $n+i-r_i$ is connected to $\bar i$, and therefore
$r_{n+i-r_i}=\infty$.  Thus a  bracket vector also satisfies condition
(ii).

Let $r$ be a vector satisfying (i) and (ii).  Recall that in  
Lemma 5, we showed that $r$ determines $C_i$ for all
$i$.  
By conditions (i) and (ii), the $C_i$ determined by $r$ do not
cross each other.  
Let $R$ be the union of the $C_i$ together with their 
symmetric partners.  The chords of $R$ divide the $(2n+2)$-gon into 
regions.  Construct a triangulation $S$ by triangulating the regions with green
chords, using the construction of  Lemma 3.  Now, for each $i$,
$C_i(S)$, if it is not an edge segment, 
is a red chord of $S$, and is therefore contained
in $R$.  Among the chords connected to $i$ and contained in $R$, 
$C_i$ is the first when encountered proceeding clockwise from $\bar 1$,
and therefore $C_i(S)=C_i$, and so $r(S)=r$, as desired.\qed  
\enddemo

We will now define an order on $\tnb$.  
For $S, T \in \tnb$, let $S\leq T$ iff for all $i$,
$r_i(S)\leq r_i(T)$.  

\proclaim{Proposition 4} The covering relations in this order on $\tnb$ 
are exactly those
described by Theorem 1. \endproclaim

\demo{Proof} We begin by proving some lemmas.

\proclaim{Lemma 6} Let $S, T\in \tnb$ such that $S \lessdot T$.  
Then there exists some $k$ such that $r_j(S)=r_j(T)$ for all $j\ne k$.
Further:

If $C_k(S)=ka$ is pure, then $C_k(T)$ connects $k$ to the endpoint of
$C_a(S)$ which is not $a$.  

If $C_k(S)=k\bar a$ is mixed, then $C_k(T)=k\bar b$ where $b$ is the
largest number $a>b>k$ such that $r_b(S)=\infty$, or $b=k$ 
if there is no such number.   

Conversely, if $S$ and $T$ are two triangulations with $r_i(S)=r_i(T)$ for
$i\ne k$ for some $k$, and $C_k(S)$ and $C_k(T)$ are related
as decribed above, then
$S\lessdot T$ in $\tnb$. 
\endproclaim

\demo{Proof}
Note that our assumption is that $S\lessdot T$ with respect to the order 
defined just above in terms of bracket vectors --- this lemma is part of the
proof that the covering relations of this partial order are as described
in the statement of Theorem 1.  
We write $s_i$ for $r_i(S)$ and $t_i$ for $r_i(T)$.

First, we prove the forward direction.  Throughout this paragraph, (i) and 
(ii) refer to conditions (i) and (ii) of Proposition 3, which describe when
an $n$-tuple is a legal $B_n$ bracket vector.  
Suppose that $s_i$ and $t_i$ coincide for $i>k$, but
$t_k>s_k$.  
We divide into cases.  
We suppose first that $C_k(S)=ka$ is pure.  
Then $t_k\geq s_k+s_{a}+1$ by
applying (i) at $(k,a)$.  
Now
$d=(s_1,\dots,s_{k-1},s_k+s_a+1,s_{k+1},\dots,s_n)$ is a valid 
 bracket vector: it satisfies (i) at 
$(k,j)$ with $k<j$ because $t$ does, while it satisfies (i) at $(j,k)$ with
$a<j<k$ because $s$ does, and it satisfies (i) at $(j,k)$ with $j<
a$ because $s$ satisfies (i) at $(j,a)$, and it clearly
satisfies (ii).  It is clear that $s <d\leq t$, so $t=d$.  The description
of $C_k(T)$ in the statement of the lemma follows immediately.

Now suppose that $C_k(S)=k\bar a$ is mixed.  
Observe first that $a>k$ because $r_k(S)<\infty$.  
Let $\bar b$ be the 
first vertex encountered counterclockwise proceeding from $\bar a$ such that
$r_b(S)=\infty$.  Set $x=n+k-b$.
If no such vertex is encountered before reaching $\bar k$, set $b=k$ and $x=\infty$.
Then let $d=(s_1,\dots, s_{k-1},x,\dots,s_n)$.
This is a valid  bracket vector.
Since $t_i=s_i$ for $i>k$, $t_k\geq x$.
Thus $t\geq d>s$, so $t=d$, as desired.    Again, the statement in the
lemma describing $C_k(T)$ follows immediately.  

Finally, we prove the converse.  Given such $S$ and $T$,
we know that $r_i(S)=r_i(T)$ for $i\ne k$, and $r_k(S)<r_k(T)$.  
We remark that there is no legal bracket vector lying betwen
$r(S)$ and $r(T)$, and we are done.  
\qed \enddemo

\proclaim{Lemma 7} Let $S\lessdot T$ in $\tnb$.  The $k$ be as in the 
statement of Lemma 6, so that $r_i(S)=r_i(T)$ for $i \ne k$.  Then:

If $C_k(S)$ is pure, then $C_i(S)=C_i(T)$ for $i\ne k$.  

If $C_k(S)$ is mixed, then $C_i(S)=C_i(T)$ for $i \ne k, i \ne b$, 
where $k\bar b=C_k(T)$.  \endproclaim

\demo{Proof} As in the previous proof, we write 
$s_i$ for $r_i(S)$ and $t_i$ for $r_i(T)$. 

Consider first the case where $C_k(S)$ is pure.  
Clearly, if $s_i=t_i\ne \infty$, then $C_i(S)=C_i(T)$.  So suppose
that $s_i=t_i=\infty$.  
In this case, recall from Lemma 5 that 
$C_i(S)=i\bar x$, where $x$ is the smallest vertex with $s_x-x\geq n-i$.
Similarly, let $C_i(T)=i\bar y$, so $y$ is the smallest vertex with
$t_y-y \geq n-i$.  We could only have that $x\ne y$ if $y=k$ and
$x\ne k$.  
Now $t_k-k=s_k+s_a+1-k=t_a-(k-(s_a+1))=t_a-a$.  Since $a<k$, $k$
cannot be the smallest vertex $y$ 
with $t_y-y\geq n-i$, so $y=k$ is impossible.  

Now consider the case where $C_k(S)$ is mixed.  Let $C_k(T)=k\bar b$.
Let $i$ be neither $k$ not $b$.  
As before, the only problematic case is when $s_i=t_i=\infty$,
and $C_i(T)=i\bar k\ne C_i(S)$.   If $i\geq a$ (and $C_i(T)=i\bar k$) then
$C_i(S)=i\bar k$ also.  
If $a<i<b$, then $s_i=\infty$ contradicts the characterization of
$b$ given in Lemma 6.  
If $i<b$, then, since $C_k(T)=k\bar b$, $k$ cannot be connected to 
$i$ in $T$, contradicting our assumption.  This completes the proof of
the lemma.  \qed \enddemo

\proclaim{Lemma 8} Let $S\lessdot T$ in $\tnb$, and let $r_j(S)=r_j(T)$ for
all $j\ne k$.  Then the only red chords of $T$ which are not red chords
of $S$ are $C_k(T)$ and its symmetric partner, and  
all the red chords of $S$ are also chords of $T$
(though not necessarily red).  
\endproclaim

\demo{Proof} 
By Lemma 7, if $C_k(S)$ is pure, then it is immediate that the only red 
chords of $T$ which are not red chords of $S$ are $C_k(T)$ and its
symmetric partner.  If
$C_k(S)$ is mixed, then if we write $C_k(T)=k\bar b$, 
it is clear that the only red chords of $T$ which are not red in $S$ are
$C_k(T)$ and $C_b(T)$ and their symmetric partners --- but in fact
$C_b(T)$ is the symmetric partner of $C_k(T)$, and the first statement
of the lemma is proved.

To prove the second statement, write $A$ for $C_k(S)$, and 
suppose first that $A=ka$ is pure.  

Let $C_k(T)$ be $kv$, where $v$ may be barred or unbarred.  
Any red chord of $S$ other than $A$ is $C_i(S)$ for some $i\ne k$, and
is therefore, by Lemma 7, also a red chord of $T$.  
We must now dispose of the chord $A$, if it is not an edge segment.
If $A$ is a red chord of $T$, we are done, so suppose otherwise.  
Consider the 
division of $T$ into regions by its red chords.  We know that $kv$ and $av$ are
red chords of $T$, so if there were a red chord of $T$ crossing $A$, it
would have to be $iv$ for some $a<i<k$.  But the fact that $s_i=t_i$ 
 would then force $C_i(S)$ to
cross $A$ also, which is impossible, since they are both chords of $S$.  

Thus, $a$ and $k$ are in the same region of $T$.  Since $k$ is the 
largest unbarred vertex in this region, $k$ and $a$ are connected by a green
chord of $T$, 
by the construction of Lemma 3.  So $A$ is a green chord in $T$.  

Suppose now that $C_k(S)=k\bar a$ is mixed.  We continue to denote
$C_k(S)$ by $A$.   
Let $C_k(T)=k\bar b$.  
Let
$D=C_b(S)$.

As before, any red chord of $S$ other than $A$ or $D$ is
also a red chord of $T$, so we need only worry about $A$ and $D$.  
Suppose $A$ is not a red chord of $T$.  Then, as in the case where 
$A$ is pure, we check that in order for a red chord of $T$ to cross
$A$, the corresponding red chord of $S$ would also cross $A$, which is 
impossible.  Since $A$ connects the largest unbarred vertex and the largest
barred vertex of the region of $T$ containing it, by the construction of Lemma
3, it is a green chord of $T$.  

Now consider $D$.  Note that $D$ is mixed, since $s_b=t_b=\infty$.  The 
argument now proceeds in the same way as for $A$.
\qed \enddemo

We now begin the proof of Proposition 4 proper.  Let $S\lessdot T$ in
the order on $\tnb$ defined by $S\leq T$ iff $r_i(S)\leq r_i(T)$ for
all $i$.  We wish to show that $S$ and $T$ are related by a diagonal flip
as in the statement of Theorem 1.  

By Lemma 8, we can consider the division of $S$ and $T$ into regions by
the red chords of $S$.  In any of these regions, the chords of $S$ are those
of the unique triangulation of the region by green chords.  Since all the
red chords of $T$ are red chords of $S$, except for $C_k(T)$ and its symmetric 
partner, the
same thing is true for $T$, except in the region(s) which contain 
$C_k(T)$ and $\overline{ C_k(T)}$.  Thus $S$ and $T$ coincide except in 
the region containing $C_k(T)$ which we denote $Z$, and the region 
containing $\overline{C_k(T)}$, which we denote $\bar Z$.  $Z$ and $\bar Z$ may
coincide.

Consider the boundary of the region $Z$.  Consider first the case where 
$A=ka$ is 
pure.  Let $C_k(T)=kv$ where $v$ may represent a barred or unbarred vertex.
Begin at $k$. Proceeding 
counter-clockwise around the boundary of $Z$, 
the next vertex is $a$, and the next is $v$.  In $S$,
all of these vertices are connected to the largest unbarred vertex of 
$Z$, say $i$.  (Note that $i\ne k$ since $C_k(S)=A$.)  In $T$, all
the unbarred vertices of $Z$ are connected to $i$ except $a$; $k$ and $v$ 
are connected by a red chord.  Thus, we see that $S$ and $T$ differ
in that the green diagonal of $ikav$ has been replaced by the red diagonal
of $ikav$,
and similarly for the symmetric partner of $ikav$.

Next consider the case where $A=k\bar a$ is mixed.  
Consider the boundary of $Z$, beginning
at $k$ and proceeding counter-clockwise.  The next vertex 
encountered is $\bar a$.  So $\bar a$ is the largest barred vertex of $Z$,
so (by the construction of Lemma 3) all the barred vertices of $Z$ are
connected to $\bar a$ in $S$.  
It follows that
none of the barred vertices of $Z$ except the smallest and the largest can
have the corresponding entries of $r(S)$ be $\infty$.  Thus, writing 
$C_k(T)=k\bar b$, we have that $\bar b$ is the
smallest barred vertex of $Z$.  Thus, $S$ and $T$ differ in $Z$ 
within the quadrilateral defined by the largest and smallest barred and
unbarred vertices of $Z$; in $S$ the larger ones are connected, while in
$T$ the two smaller ones are conected.  If $Z \ne \bar Z$, then the same
analysis holds in $\bar Z$.

Thus, we have shown that if $S\lessdot T$ in $\tnb$, then they are related
by a minimal flip as in Theorem 1.  We must now check that if $S$ and
$T$ are related by a minimal flip as in Theorem 1, then $S\lessdot T$ in
$\tnb$.  

So suppose that $S$ and $T$ are related by a minimal flip as in Theorem 1:
that is to say, there is a chord $C$ of $S$ which is green, such that 
$T$ can be obtained from $S$ by replacing $C$ and $\bar C$ by the
other diagonals of $Q(C)$ and $Q(\bar C)$.  
It is a case-by-case check, based on the positions of the four vertices of
$Q(C)$, that $r_i(S)=r_i(T)$ for $i\ne k$ for some $k$.  One then checks
that $C_k(S)$ and $C_k(T)$ are related as in Lemma 6.  By the converse 
direction of Lemma 6, this then implies that $S\lessdot T$ in $\tnb$.  
This completes the proof of Proposition 4.  \qed \enddemo

Our next goal is to prove that the $B_n$ Tamari order is really a lattice.
Before we can prove that, we need some preliminary results.  

Let $M_n$ denote the $n$-tuples with entries in $[0,n-1]\cup \{\infty\}$,
with the Cartesian product order.    
Let $\mi$ denote the elements of $M_n$ which satisfy 
condition (i) of Proposition 3.  Let $\mii$ denote the elements of $M_n$ which satisfy
condition (ii) of Proposition 3.

\proclaim{Proposition 5} There exist maps $\upa\,: \mii \rightarrow
\tnb$, $\downa\,: \mi \rightarrow \tnb$, which satisfy  the following
conditions:
$$\align
&f\leq r(S) \text{ iff } \upa(f)\leq S \tag1\\
&r(S)\leq f \text{ iff } S\leq \downa(f).\tag2\endalign $$ 
\endproclaim

\demo{Proof} 
Let $f\in \mii$.  Define $g\in M_n$  inductively, as follows:
$$g_i=\max_{j\leq f_i} (g_{i-j}+j).$$
By construction, $g$ satisfies (i), and $g$ satisfies (ii) since $f$
does.  Thus, $g$ is a bracket vector.  Let $\upa(f)$ be the corresponding
triangulation. 
Now statement (1) is clear.

Let $f\in \mi$.  Define $g\in M_n$, as follows: set $g_i=f_i$, unless
$f_i\geq i$, and $f_{n+i-f_i}\ne \infty$.  In this case, set $g_i$ to be the 
largest number less than $f_i$ such that $f_{n+i-g_i}=\infty$ or $g_i<i$.  
By construction, $g$ satisfies (ii), and it is a straightforward check 
that $g$ will also satisfy (i), since $f$ does.  (2) is also clear.\qed 
\enddemo

Using these maps, we can prove that meet and join exist in 
$\tnb$ by giving simple descriptions of them.  

\proclaim{Proposition 6} The Tamari order on $\tnb$ is a lattice.  
The lattice operations are as follows:
For $S,T \in \tnb$, $S\vee T = \upa\!(\max(r(S),r(T)))$ and  $S\wedge
T=\downa(\min(r(S),r(T)))$.    
\endproclaim

\demo{Proof} It is clear that, in $M_n$, the join of $r(S)$ and $r(T)$ 
is $\max(r(S),r(T))$.  Now,
since $\max(r(S),r(T))\in \mii$, for any $W\in \tnb$, $W\geq S$ and $W\geq T$
iff $r(W)>\max(r(S),r(T))$ iff $W>\upa(\max(r(S),r(T)))$, so 
$\upa(\max(r(S),r(T)))=S\vee T$.  The same argument holds for 
$S\wedge T$, once we observe that $\min(r(S),r(T))\in \mi$.  \qed\enddemo

This completes the proof of Theorem 1.  The Hasse diagram of $T^3_B$ is
shown in Figure 5, at the end of the paper.  \qed 

\head Noncrossing partitions \endhead

The $A_n$ noncrossing partitions, $\ncna$, are partitions of $n+1$ into 
sets such that if $v_1,\dots,v_{n+1}$ are $n+1$ points on a circle,
labelled in cyclic order, and if $B_1,\dots, B_r$ are the convex hulls
of the sets of vertices corresponding to the 
blocks of the partition, then the $B_i$ are non-intersecting. 

There is a bijection from $\tna$ to 
$\ncna$
as follows.  For $S \in \tna$, erase all the green chords and exterior 
edges of $S$, together with the
vertices $0$ and $n+2$.
Then move
the endpoints of each red chord $ij$ 
a little bit, the lower end point a little
clockwise,
the higher endpoint a little counterclockwise (so $i$ and $j$ are both on the
upper side of the chord).  These chords now divide the vertices in $[n+1]$
into subsets, which form a noncrossing partition by construction.  
Figure 3 shows the triangulation from Figure 1, together with the 
noncrossing partition which it induces: $\{14,23,5\}$.

$$\epsfbox{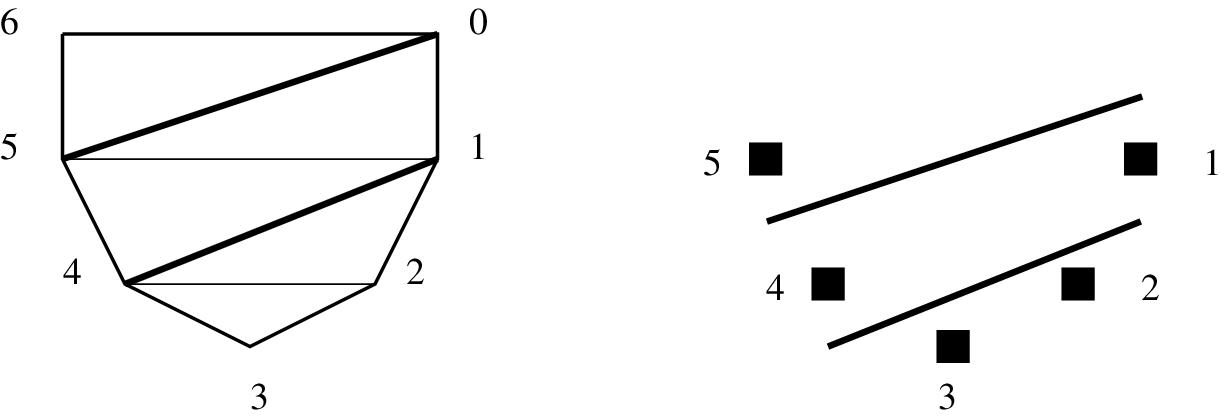}$$
$$\text{Figure 3}$$

Note that the noncrossing partitions are often considered as being
ordered by refinement; this order is quite different from the Tamari
order. 

As defined by Reiner [Rei], 
the $B_n$ noncrossing partitions, $\ncbn$, are partitions of the set
$1,\dots,n,$ $\bar 1,\dots, \bar n$, which have the properties that the 
partition remains fixed under interchanging barred and unbarred elements,
and that if $2n$ points are chosen around a circle and labelled 
cyclically $v_1,\dots,v_n,v_{\bar 1},\dots,v_{\bar n}$, then the 
convex hulls of the vertices corresponding to the blocks of the partition do not intersect.  

We now define a map 
$\psi$ from $\tnb$ to $\ncbn$, analogous to 
that in type $A$.  
Erase all green chords and exterior edges.  
Move both endpoints of mixed red chords 
slightly counterclockwise.  
Move the 
endpoints of pure red chords slightly together 
(so that the vertices both lie on the side
of the chord which includes the larger part of the polygon).  
Erase the vertices $n+1$ and $\overline{n+1}$.
The remaining vertices are now partitioned by the red chords, in what is
clearly a $B_n$ noncrossing partition.  Figure 4 shows the triangulation
from Figure 2, together with the $B_n$ noncrossing partition which it
induces: $\{ 1\bar 2 \bar 5 \bar 6, 34, \bar 1 2 5 6, \bar 3 \bar 4\}$.

$$\epsfbox{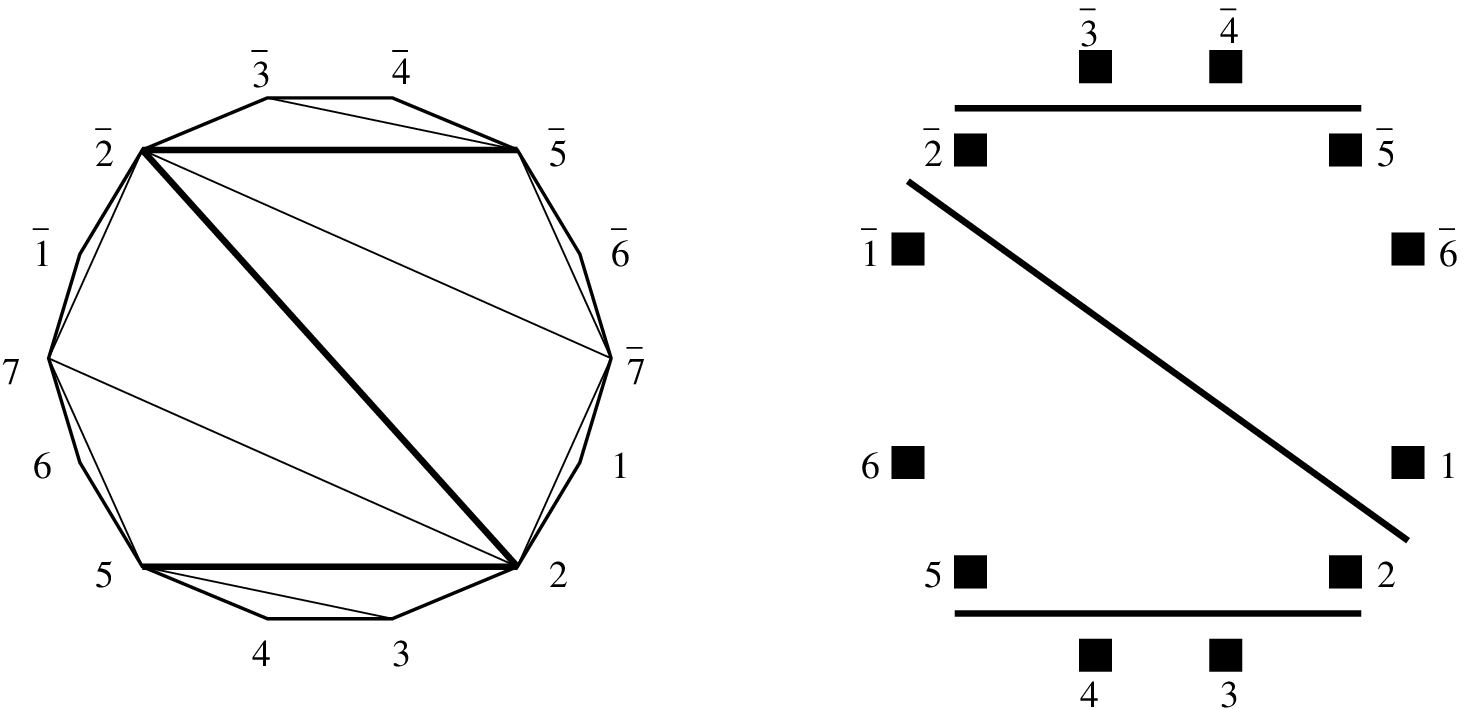}$$
$$\text{Figure 4}$$

\proclaim{Proposition 7}
The map $\psi$ is a bijection from 
$\tnb$ to $\ncbn$. \endproclaim

\demo{Proof} We remark first that $\tnb$ and $\ncbn$ have the same 
cardinality, $\binom {2n}{n}$ (see [Sim] and [Rei] respectively).  
Thus, it suffices to show that $\psi$ is an injection.

Fix $S\in\tnb$.  Let $T$ be a triangulation in the same fibre as $S$.  
Pick $i \in [n]$.  Starting at $i-1$, search counter-clockwise around the 
$(2n+2)$-gon for the first vertex whose label is in the same block as $i$
in $\psi(S)$.  Let this vertex be $v$.  

Consider first the case that $v=i-1$.  Then there must be no red chords 
connected to $i$ in $T$, so $C_i(T)$ is an edge segment.  
Next consider the case 
where $i-1>v\geq 1$.  Then there must be a red chord $iv$ in $T$.  This 
chord cannot be $C_v(T)$, so it must be $C_i(T)$.  Finally consider the case
where $v=\bar j$.  Then $T$ must contain the red chord 
$i\left(\overline{j+1}\right)$, 
and it cannot contain any chord $i\bar k$ with $k<
j+1$.  Thus $C_i(T)$ must be $i\left(\overline{j+1}\right)$.  
Similarly, if $v=n$, then $C_i(T)$ must be $i\bar 1$.  
Finally, if $v$ is unbarred and $n>v\geq i$, then there must be a red chord
$(i-1)(v+1)$, so $C_i(T)$ is an edge segment.  

Since the $C_i(T)$ suffice to determine $T$, it follows that $S$ and $T$
coincide, and $\psi$ is an injection, as desired.  \qed\enddemo

\head EL-Shellability \endhead

Recall that an element $x$ of a lattice $L$ is said to be left modular if, 
for all $y<z \in L$,
$$(y \vee x)\wedge z=y\vee(x\wedge z).\tag3$$  
In this section we shall
prove the following theorem:

\proclaim{Theorem 2} $\tnb$ has an unrefinable chain of left modular elements.  \endproclaim

The analogous fact that 
$\tna$ posesses a unrefinable chain of left modular elements
was first proved by Blass and Sagan [BS].  

It was shown in [Liu] that a
lattice having an unrefinable chain of left-modular elements has an
EL-labelling.  In particular,
this shows that 
the order complex of any interval in such
a lattice is shellable and hence contractible or 
 homotopy equivalent to a wedge of spheres. For more on 
EL-labelling and EL-shellability, see [Bj\"o].  Neither 
$\tna$ nor $\tnb$ are graded, and therefore the EL-shellings are not pure.
For more on ungraded EL-shellability, see [BW1, BW2].  

Thus, Theorem 2 implies the following corollary:

\proclaim{Corollary} $\tnb$ is EL-shellable.\endproclaim

\demo{Proof of Theorem 2} 
For $1\leq i \leq n$ and $t\in [1,n-1] \cup \{\infty\}$, let  
$S_{i,t}$ denote the triangulation with  bracket vector as follows:
$$r_j(S_{i,t})= \left\{ \matrix  0& \text{ for } j<i\\
t & \text{ for } j=i\\
\infty & \text{ for } j>i. \endmatrix \right.$$

\proclaim{Lemma 9} 
$S_{i,t} \in \tnb$ is left 
modular.\endproclaim

\demo{Proof} The proof is just a calculation, verifying (3).
Let $Y<Z \in \tnb$.  Let $r(Y)=(y_1,\dots,y_n)$, 
$r(Z)=(z_1,\dots,z_n)$.  For ease of notation, we split into three cases,
when $t=\infty$, when $t\leq i-2$, and when $i-1\leq t<n$.  

Suppose first that $t=\infty$.  We observe that 
$$r(Z\wedge(S_{i,\infty}\vee Y))=(y_1,\dots,y_{i-1},z_i,\dots,z_{n})=
r((Z\wedge S_{i,\infty})\vee Y),$$
as desired.  

Suppose next that $t \leq i-2$.  Applying Proposition 5, we see that:
$$ \align
r(S_{i,t} \vee Y)&= (y_1,\dots,y_{i-1},
\max_{0\leq j\leq t} y_{i-j}+j,\infty,\dots,\infty)\\
r(Z\wedge(S_{i,t}\vee Y))&=(y_1,\dots,y_{i-1},
\min(z_i,\max_{0\leq j\leq t} (y_{i-j}+j)),z_{i+1},\dots,z_n)\\
r(Z \wedge S_{i,t})&=(0,\dots,0,\min(z_i,t),z_{i+1},\dots,z_n)\\
r((Z \wedge S_{i,t})\vee Y)
&=(y_1,\dots,y_{i-1},\max_{0\leq j \leq \min(z_i,t)} (y_{i-j}+j),
z_{i+1}.\dots,z_n).\endalign$$

Thus, we need only check that
$$\min(z_i,\max_{0\leq j\leq t} (y_{i-j}+j))=
\max_{0\leq j \leq \min(z_i,t)} (y_{i-j}+j).$$

For $j\leq z_i$, we have that 
$z_i\geq z_{i-j} +j\geq y_{i-j}+j$,  
so 
$$\max_{0\leq j \leq \min(z_i,t)} (y_{i-j}+j)\leq z_i.$$
Thus, we alter nothing by rewriting:
$$\max_{0\leq j \leq \min(z_i,t)} (y_{i-j}+j)=
\min(z_i,\max_{0\leq j \leq \min(z_i,t)} (y_{i-j}+j)).$$
If $j>z_i$, then $y_{i-j}+j>z_i$, so 
$$\min(z_i,\max_{0\leq j \leq \min(z_i,t)} (y_{i-j}+j))=
\min(z_i,\max_{0\leq j \leq t} (y_{i-j}+j))$$
and we are done.  

The case where $i-1 \leq t <n$ is similar.  Let $t'$ denote the largest integer
less than $t$ and greater than $i-1$ such that $z_{n+i-t'}=\infty$.  
If there is none, set $t'=i-1$.  
Then
$$\align
r(S_{i,t}\vee Y)&=(y_1,\dots,y_i,\max(t,\max_{0\leq j< i} (y_{i-j}+j)),
\infty,\dots,\infty) \\
r(Z\wedge(S_{i,t}\vee Y))&=(y_1,\dots,y_i,\min(z_i,\max(t',
\max_{0\leq j< i} (y_{i-j}+j))),z_{i+1},z_n)\endalign$$
where $t$ has turned into $t'$ because of the application of $\downa\,$.

On the other hand
$$\align
r(Z\wedge S_{i,t})&=(0,\dots,0,\min(z_i,t'),z_{i+1},\dots,z_n)\\
r((Z\wedge S_{i,t})\vee Y)&=(y_1,\dots,y_{i-1},\max(\min(z_i,t'),\max
_{0\leq j \leq \min(i-1,z_j)} (y_{i-j}+j)),\\
&\qquad\qquad\qquad z_{i+1},\dots,z_n).\endalign$$

The remainder of the argument is similar to the previous case.  
This completes the proof of Lemma 9.\qed\enddemo

Now observe that $\hat 0 \lessdot S_{n,1}\lessdot S_{n,2}
\lessdot\dots\lessdot S_{n,\infty}\lessdot S_{n-1,1}\lessdot \dots
\lessdot S_{1,\infty}=\hat 1$ forms an unrefinable chain in $\tnb$.  This proves 
Theorem 2.  \qed\enddemo

\head Homotopy types of intervals \endhead

As we have already remarked, the fact that $\tnb$ is EL-shellable 
implies
that the order complex of any interval is either contractible or has
the homotopy type of a wedge of spheres.  In this section, we shall 
prove that it is in fact either contractible or homotopic to a 
single sphere.  One reason that such a result is of interest is that
it implies that the M\"obius function of any interval in $\tnb$
is $0$, $-1$, or $1$.  

\proclaim{Theorem 3} The order complex of an 
interval in $\tnb$ is either contractible or homotopy equivalent to a 
single sphere.\endproclaim

\demo{Proof} An element of a lattice called a join irreducible if it 
cannot be written as the join of two stricty smaller elements.  
We now prove
some results concerning the
join irreducibles of $\tnb$.  

For $1\leq t\leq i-1$, let $W_{i,t}$ denote the triangulation 
whose  bracket vector consists of $t$ in the $i$-th place, all the
other entries being zero. 

For $i\leq t<n$, let $W_{i,t}$ denote the triangulation defined by:
$$  r_j(W_{i,t})=\left\{ \matrix t &\text { for }  j=i\\
                                \infty &\text { for }  j=n+i-t \\
                                0 & \text { otherwise} \endmatrix
\right.$$

Let $W_{i,\infty}$ denote the triangulation whose bracket vector 
consists of a single $\infty$ in the $i$-th place, all the other entries
being zero.         

Write $\Cal{W}$ for the set of all the $W_{i,t}$.  
     
\proclaim{Proposition 8} The join irreducibles of $\tnb$ are 
exactly $\Cal{W}$. 
\endproclaim

\demo{Proof} It is easy to see that these elements are join irreducible
and that any element of $\tnb$ can be written as the join of those
$W_{i,t}$ below it.\qed \enddemo

For $S<T$, write $\Cal{W}(S,T)$ for the set of join irreducibles of 
$T^B_n$ which lie below $T$ but not below $S$.

We now define the EL-labelling 
$\gamma$ of [Liu].  Let $L$ be a 
lattice, and let $\hat 0=x_0\lessdot x_1 \lessdot \dots \lessdot x_r=\hat 1$
be an unrefinable chain of left modular elements.  Let $\Cal{W}_i$ be the 
set of join irreducibles below $x_i$ but not below $x_{i-1}$.  For 
$y<z$ in $L$, let $\W(y,z)$ be the set of irreducibles below $z$ but not
below $y$.  For any $S\lessdot T$ in $L$, label the corresponding edge 
of the Hasse diagram by:
$$\gamma(S,T)=\min\{i\mid \W_i \cap \W(S,T)\ne \emptyset\}.$$

\proclaim{Proposition 9 [Liu]} For $L$ a lattice with an unrefinable left 
modular chain, the labelling $\gamma$ defined above is an EL-labelling.
\qed \endproclaim

We now interpret this labelling in our context.  
Observe that $W_{i,t}$ is that unique join irreducible that
lies below $S_{i,t}$ but not below those $S_{j,p}$ below
$S_{i,t}$.  Instead of numbering the $S_{i,t}$, we proceed as follows.
  We put a new linear order, $\prec$, on 
$\Cal{W}$, so that:
$$W_{n,1}\prec W_{n,2}\prec \dots\prec W_{n,\infty}\prec W_{n-1,1}\prec\dots
\prec W_{1,\infty},$$
in other words, so that:
$$W_{i,t}\prec W_{j,p} \text{ iff } S_{i,t}<S_{j,p}.\tag4$$
Now, we label the edges of the Hasse diagram of $\tnb$ by join irreducibles: 
if $S\lessdot T$, we label the edge $(S,T)$ by 
the minimal element of $\W(S,T)$ (under $\prec$).  Clearly, this
is equivalent to the labelling defined by [Liu], and is therefore
an EL-labelling.  

\proclaim{Lemma 10} If $Y\lessdot Z$, such that $r_i(Y)=r_i(Z)$ for
$i\ne k$, then $\gamma(Y,Z)=W_{k,t}$ for some $1\leq t\leq \infty$. 
\endproclaim

\demo{Proof} Consider separately the cases where $r_k(Y)\leq k-2$, where
$k-1\leq r_k(Y)\leq \infty$, and where $r_k(Y)=\infty$.\qed \enddemo

Recall from [BW2]
that given a poset with an EL-labelling, the order 
complex of an interval $[y,z]$ is homotopic to a wedge of spheres,
one for each unrefinable chain from $y$ to $z$ such that the labels
strictly decrease as one reads up the chain.  Such chains are called
{\it decreasing chains}.

Thus, Theorem 3 will follow from the following lemma:

\proclaim{Lemma 11} For $Y<Z \in \tnb$, there is at most one 
decreasing chain from $Y$ to $Z$. \endproclaim

\demo{Proof}
Let $Y=T_0\lessdot T_1\lessdot \dots \lessdot T_r=
Z$ be a decreasing chain from $Y$ to $Z$.   

For each $i$, let $k_i$ denote the unique place (provided by 
Lemma 6) where $r(T_i)$ and $r(T_{i+1})$ differ.  By Lemma 11,
the label on the edge $(T_i,T_{i+1})$ is $W_{k_i,t_i}$ for some $t_i$.  
Since the labels are decreasing by assumption, $k_0$ must be the 
first index where $r(Y)$ and $r(Z)$ differ.  Thus $r(T_1)$ must be the
smallest legal bracket vector obtainable by increasing the $k_0$ position
of $r(T_0)$.  This determines $T_1$ uniquely, and the remaining 
$T_i$ are determined inductively, proving the lemma.
\qed\enddemo

This completes the proof of Theorem 3. \qed\enddemo

An exact description of when the homotopy type of the interval $[Y,Z]$
in $\tnb$ is homotopic to a sphere, and when it is contractible
(expressed in terms of the bracket vectors of $Y$ and $Z$), 
has been worked out in
[Sa2].

\head Generalizing to Type $BD_n^S$ \endhead

Here we fix $n$ and a subset $S$ of $[n]$.  We will be operating in
type $BD_n^S$, a notation introduced in [Rei] which we now explain.
This is not a type in the usual sense.  Rather, it refers to
a certain
hyperplane arrangement between those associated to $B_n$ and $D_n$.

Recall that a root system gives rise to a 
hyperplane arrangement by taking all the hyperplanes through the origin
perpendicular to 
roots.  The $B_n$ arrangement 
therefore consists of all those
hyperplanes defined by $x_i\pm x_j=0$, together with those
defined by $x_i=0$, for $1\leq i,j \leq n$, 
while the $D_n$ arrangement consists only of those hyperplanes defined 
by $x_i\pm x_j=0$ for $1\leq i,j\leq n$.  
Now, for $S \subset [n]$, the
$BD_n^S$ hyperplane arrangement consists of those hyperplanes
defined by $x_i\pm x_j=0$ together with $x_i=0$ for $i \not\in S$.  
When $S = \emptyset$ we recover the $B_n$ arrangement, 
while if $S=[n]$ we obtain
the $D_n$ arrangement.

The $B_n$ partitions, $\Pi_n^B$, 
are by definition those partitions of the 
set $\{1,\dots,n,$ $\bar 1,\dots,\bar n\}$ which are fixed under the map
interchanging $i$ and $\bar i$, and such that there is at most one block
which contains any $i$ and $\bar i$ simultaneously.  
This is a suitable definition of $\Pi_n^B$ because its elements are 
naturally in bijection with the elements of the intersection lattice of
the $B_n$ arrangement. 
$\ncnb$ is a subset
of $\Pi_n^B$. 

The intersection lattice of the 
$BD_n^S$ hyperplane arrangement is a subset of that of type $B_n$.  This
allows a natural definition of $BD_n^S$ partitions, $\Pi_n^S$, as a subset of
$\Pi_n^B$.  By this approach, one obtains that 
$\Pi_n^S$ consists of those partitions of $\Pi_n^B$ which
do not contain any block consisting solely of $\{i,\bar i\}$ for $i \in S$.  
In [Rei], Reiner defined $\ncbd$, the noncrossing partitions of 
type $BD_n^S$, 
by $\ncbd=\ncbn\cap\Pi_n^S$.  

In particular, Reiner defined the noncrossing
partitions of type $D_n$ to be the noncrossing partitions of type
$BD_n^{[n]}$.  
Since [Rei] was written, it has become clear that there is a more natural 
definition of noncrossing partitions of type $D$, see Athanasiadis and Reiner
[AR].  Thus, even though type $BD_n^{[n]}$ is associated to the $D_n$
hyperplane arrangement, it should not be confused with type $D_n$.

Let $\tns$ be those triangulations which correspond under
$\psi$ to partitions in $\ncns$.  
We can describe them more directly as follows:

\proclaim{Lemma 12} $\tns$ consists of those triangulations which do not 
contain the triangles $i, i+1, \overline{i+1}$ and $\bar i, i+1, 
\overline{i+1}$ for any 
$i \in S$.  $\tns$ can also be characterized as the set of triangulations
$T$ such that
$r_i(T) \ne n-1$ for any $i\in S$.  
\endproclaim

\demo{Proof} We prove the first statement by showing 
that $\psi(T)$ contains the block $\{i,\bar i\}$ iff
$T$ contains the triangles $i, i+1, \overline{i+1}$ and $\bar i, i+1, 
\overline{i+1}$.  

If $\psi(T)$ contains the block $\{i,\overline i\}$, then
$T$ contains red chords $\bar i (i+1)$ and 
$i \left(\overline{i+1}\right)$.  Therefore, $T$ must also contain
either $i\bar i$ or $(i+1) \left(\overline{i+1}\right)$.  The first of 
these chords would be red, which contradicts the presence of $\{i,\bar i\}$
in $\psi(T)$.  So $T$ must contain also contain 
the edge $(i+1) \left(\overline{i+1}\right)$, so it contains the two desired
triangles.  

Conversely, suppose $T$ contains the two desired triangles.  It is immediate
that $(i+1)\left(\overline{i+1}\right)$ will be green while 
$i \left(\overline{i+1}\right)$ and $\bar i (i+1)$ will be red, and therefore
$\psi(T)$ contains the block $\{i,\bar i\}$.  This completes the proof of
the first statement.  

The second statement follows immediately from the first.  \qed\enddemo

The remainder of the paper is devoted to the proof of the following 
theorem, which generalizes Theorems 1, 2, and 3 to the broader context 
of type $BD_n^S$.
 
\proclaim{Theorem 4} $\tns$ admits a lattice structure which is a quotient 
of that on $\tnb$.  $\tns$ posesses 
an unrefinable chain of left modular elements, which implies that it is
EL-shellable.  Further, the order complex of any interval is either 
contractible or homotopic to a single sphere.  \endproclaim

\demo{Proof}
We define an equivalence relation $\ltilde_S$ on $\tnb$ as follows: 
two non-identical triangulations are equivalent iff 
they differ in that one of them, say $T$,
is not in $\tns$, and the other is the triangulation obtained by removing the
diameter of $T$ and replacing it with the other possible diameter.  It is immediate from Lemma 12 that this triangulation 
will be in $\tns$.  

We can express equivalence in terms of bracket
vectors by saying that $V$ and $W$ are equivalent if there is some $k\in S$
such that
$r_j(V)=r_j(W)$ for all
$j\ne k$, and $r_k(V)=n-1$ while $r_k(W)=\infty$.

An equivalence relation $\ltilde$ on a lattice $L$ is said to be a congruence relation
if the lattice operations pass to equivalence classes.  In this case,
 there is an induced lattice structure on the equivalence classes 
(see [Gr\"a]).  

\proclaim{Lemma 13} The relation $\ltilde_S$ on $\tnb$ is a congruence
relation. \endproclaim

\demo{Proof} It suffices to show that for $V\ltilde_S W \in \tnb$ and 
$Z\in \tnb$, $V\vee Z \ltilde_S W\vee Z$ and $V\wedge Z \ltilde_S 
W \wedge Z$.  But these are both clear from the descriptions of the
lattice operations in $\tnb$ in terms of  bracket vectors.\qed\enddemo 

Since the equivalence classes of $\ltilde_S$ each contain a single element
of $\tns$, the induced lattice structure on $\tnb/\ltilde_S$ gives rise
to a lattice structure on $\tns$.  Let us write $<_S$ for the order induced
in this way on $\tns$.

There is another way to induce a poset structure on $\tns$, namely 
that induced by its inclusion in $\tnb$.  Let us write $<_B$ for
the poset structure induced on $\tns$ in this way.  
As we now prove, $<_B$ and $<_S$  coincide (and so, once we have 
proved the lemma, we will drop the subscripts).

\proclaim{Lemma 14} The two poset structures on $\tns$, $<_B$ and $<_S$,
coincide.  \endproclaim

\demo{Proof} Observe that $\tns$ is closed under the meet in $\tnb$.  
Thus, $(\tns,<_B)$ has a meet, namely
the meet in $\tnb$  restricted to $\tns$.  
But this coincides with the 
meet defined for $(\tns,<_S)$.  Knowing that meets exist in both
poset structures, and that they coincide, implies that the two 
structures themselves coincide. \qed\enddemo

We shall generally prefer to consider the poset structure on $\tns$ as 
being induced by $\ltilde_S$.  Note that $\tns$ does not form a sublattice
of $\tnb$, because $\tns$ is not closed under the join in $\tnb$.    

It is immediate that the property of being left modular passes to equivalence
classes, so $\tns$ has a unrefinable chain of left modular elements,
and is therefore EL-shellable.  This unrefinable chain is shorter than that
of $\tnb$, because $S_{i,n-1}\ltilde_S S_{i,\infty}$ for $i\in S$.  

It is easy to see that the join irreducibles of $\tns$ are those $W_{i,t}$
such that either $i \not \in S$ or $t \ne n-1$; again, they are in bijection
with the elements of the left modular chain.  Write $\Cal{W}^S$ for the set of
all the join irreducibles of $\tns$.  As before, we define an  order 
$\prec $ on $\Cal{W}^S$ so that (4) is satisfied.  
In fact, this is just the order
induced on $\Cal{W}^S$ from its inclusion in $(\Cal{W},\prec)$.  

\proclaim{Lemma 15} If $Y \lessdot Z$ in $\tns$, then either 
$Y \lessdot Z$ in $\tnb$ or there exists some $Z'$ such that
$Y \lessdot Z' \lessdot Z$ in $\tnb$, and $Z' \ltilde_S Z$.
In the latter case, the bracket vectors of $Y$, $Z'$, and $Z$ all
coincide except in one place, where $Z'$ has $n-1$ and $Z$ has $\infty$.
\endproclaim

\demo{Proof} Suppose that $Y \lessdot Z$ in $\tns$ but 
$Y \not \lessdot Z$ in $\tnb$.  So there exists
some $Z'$ in $\tnb$, $Y<Z'<Z$.  When we pass to equivalence classes with
respect to $\ltilde_S$, one of the strict inequalities becomes an
equality, so either $Z'\ltilde_S Y$ or 
$Z' \ltilde_S Z$.  The first case is impossible, because an element of 
$\tns$ is the top element of its equivalence class mod $\ltilde_S$.  
Thus, the second case holds, which proves that $Z'\ltilde_S Z$, and therefore
that $Z' \lessdot Z$, and, since we have shown that
the only element in $\tnb$ strictly between $Y$ and $Z$ is $Z'$,
it follows that $Y\lessdot Z'\lessdot Z$.  

The bracket vectors of 
$Z'$ and $Z$ differ in one place, say $k$, and clearly 
$r_k(Z')=n-1$ and $r_k(Z)=\infty$.  The bracket vector of
$Y$ differs from that of $Z'$ in one place, and since $Y \in T^S_n$, 
$r(Y)_k \ne n-1$, so $Y$ differs from $Z'$ in the $k$-th place.  
\qed\enddemo

For $Y<Z \in \tns$, 
write $\W^S(Y,Z)$ for the set of join irreducibles of $\tns$ below $Z$ but
not below $Y$.  
Now, for $Y \lessdot
Z$ in $\tns$, we label the edge $(Y,Z)$ by the first element (with
respect to $\prec$) of $\W^S(Y,Z)$.  By the result of Liu already cited, this
is an EL-labelling.  

Using Lemma 15, we see that 
Lemma 10 holds in $T^S_n$; the same proof goes through.  The argument used to
prove Lemma 11
now goes through to show that there can be at most one decreasing chain from
$Y$ to $Z$ in $T^S_n$.  The final statement of Theorem 4 follows.  
\qed\enddemo

\head Acknowledgements \endhead

I would like to thank Nathan Reading for suggesting
the possibility of a type $B$ Tamari lattice to me, and for his amicable 
approach to the overlap in our investigations.  
I would also like to thank Vic Reiner and Christos Athanasiadis for 
discussions in which many of the ideas that appear in this paper were 
formed.  I would like to thank Marcelo Aguiar,  Nathan
Reading, and Vic Reiner for their comments on a previous version of
the manuscript.  Nirit Sandman and a pair of anonymous referees deserve
particular thanks for their close reading of a previous version of the
manuscript, which uncovered certain points requiring improvement.

\newpage

$$\epsfbox{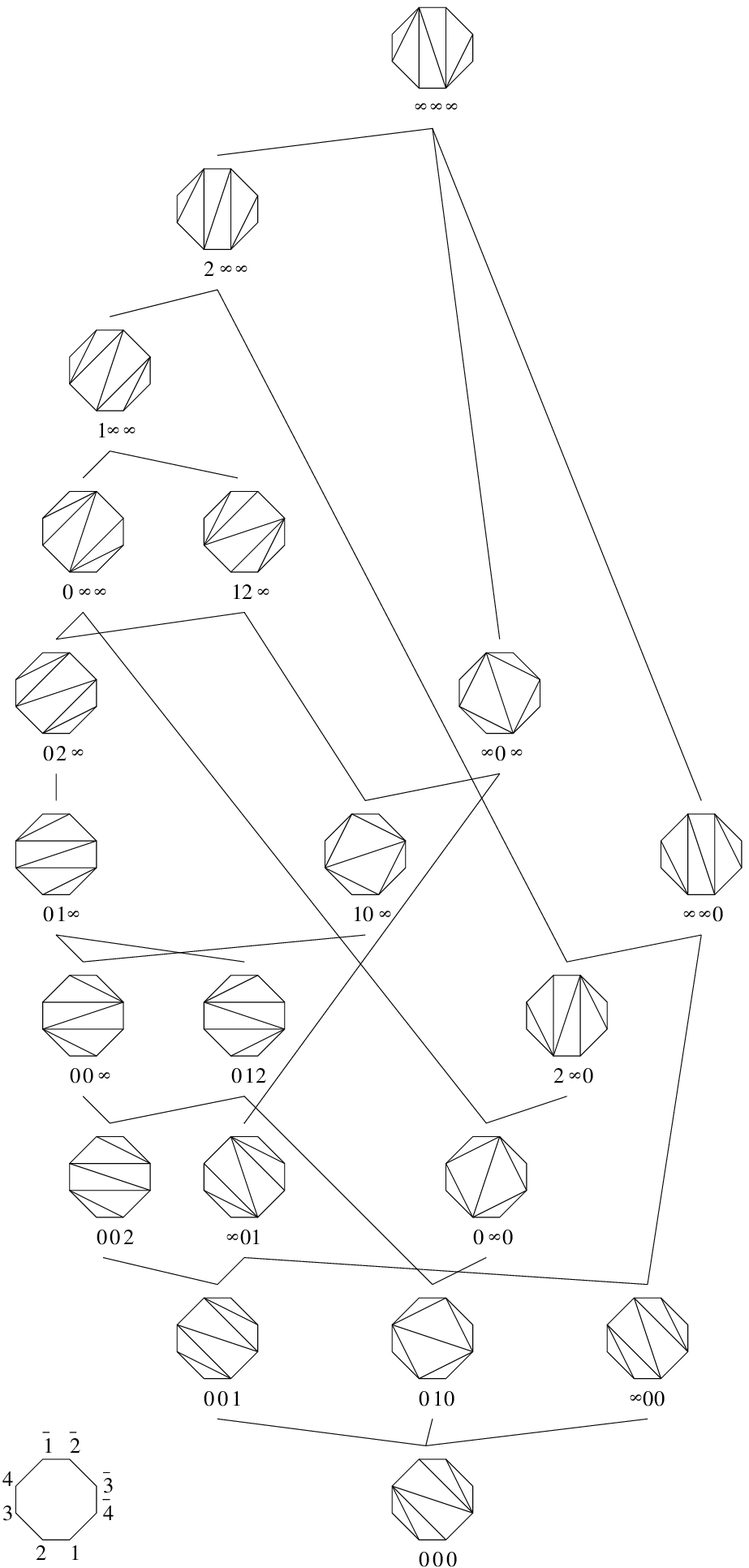}$$
$$\text{Figure 5}$$

\enddocument